\documentclass[a4paper,12pt]{article}
\usepackage{amssymb,latexsym,amsmath,amsthm,amscd,graphpap}
\usepackage{amsmath}
\usepackage{amsfonts}
\usepackage{amssymb}
\usepackage[dvips]{epsfig}
\usepackage{subfigure,array}

\newtheorem{theorem}{Theorem}[section]
\newtheorem{prop}[theorem]{Proposition}
\newtheorem{lemma}[theorem]{Lemma}
\newtheorem{coro}[theorem]{Corollary}
\newenvironment{demo}{ \noindent
\emph{\textbf{Proof~:}}}{\hfill$\square$\\
\vspace{0.4cm}}

\newcommand{\RR}{\mathbb{R}}

\newcommand{\NN}{\mathbb{N}}

\newcommand{\Cc}{\mathcal{C}}

\newcommand{\Oc}{\mathcal{O}}

\newcommand{\eps}{\varepsilon}

\newcommand{\supp}{\mathrm{supp}}
\newcommand{\no}{n$^{\text{o}}$}
\newcommand{\grad}{\nabla}

\numberwithin{equation}{section}

\newdimen\texpscorrection
\newdimen\figcenter
\def\figurewithtex #1 #2 #3 #4 #5\cr{\null
  {\goodbreak\figcenter=\hsize\relax
  \advance\figcenter by -#4truecm
  \divide\figcenter by 2
  \begin{figure}[hbt]
  \vskip #3truecm\noindent\hskip\figcenter
  \includegraphics{#1}{\hskip\texpscorrection\input #2 }
  \vskip 0.8truecm{\baselineskip=0.8\baselineskip
  \noindent \vbox{\noindent {\footnotesize #5}}\par}
  \end{figure}}}
\def\point#1 #2 #3 {\rlap{\kern #1 truecm
\raise #2 truecm \hbox{#3}}}

\begin{document}

\title{{\bf\huge How opening a hole affects the sound of a flute\\}{\Large \bf A
one-dimensional mathematical model for a tube with a small hole pierced on its
side}}

\author{Romain Joly \\ {\small Institut Fourier}\\  {\small UMR 5582
CNRS/Universit\'e de Grenoble}\\ {\small 100, rue des maths B.P. 74}\\ {\small
38402
Saint-Martin-d'H\`eres, France}\\ {\small \tt Romain.Joly@ujf-grenoble.fr}\\}

\date{May 2011}

\maketitle

\vspace{1cm}

\begin{abstract}
In this paper, we consider an open tube of diameter $\varepsilon>0$, on the
side of which a small hole of size $\varepsilon^2$ is pierced. The
resonances of this tube correspond to the eigenvalues of the Laplacian operator
with homogeneous Neumann condition on the inner surface of the tube and 
Dirichlet one the open parts of the tube. We show that this
spectrum converges when $\varepsilon$ goes to $0$ to the spectrum of an
explicit 
one-dimensional operator. At a first order of approximation, the limit spectrum
describes the note produced by a flute, for which one of its holes is open.
\\[3mm]
\end{abstract}


\section{Introduction and main result}

In this paper, we obtain a one-dimensional model for the resonances of a tube
with a small hole pierced on its side. Our arguments are based on recent thin
domain techniques of \cite{Murat}. We show that this kind of techniques applies
to
the mathematical modelling of music instruments.

\vspace{3mm}

{\noindent \bf Basic facts on wind instruments.}\\
The acoustic of flutes is a large subject of research for acousticians.
Basically, a flute is the combination of an exciter which creates a periodic
motion (a fipple, a reed etc.) and a tube, whose first mode of resonance
selects the note produced. Studying the acoustic of a flute combines a
lot of problems as the influence of the shape of the tube, the study of the
creation of oscillations by blowing in the fipple\ldots see \cite{Rossing},
\cite{Fletcher}, \cite{Dickens}, \cite{Wolfe-web} and \cite{Bolton} for nice
introductions. In this paper, we will
not consider the creation of the periodic excitation, we rather want to study
mathematically the resonances of the tube of the flute and how an open hole
affects it. Therefore, we simplify the problem by making the following usual
assumptions:\\
- the pressure of the air in the tube follows the wave equation and therefore
the resonances of the tube are the squareroots of the eigenvalues of the
corresponding
Laplacian operator.\\
- on the inner surface of the tube, the pressure satisfies homogeneous Neumann
boundary condition.\\
- where the tube is open to the exterior, we assume that the pressure is equal
to the exterior pressure which may be assumed to be zero without loss of
generality.\\

\begin{figure}[tp]
\begin{center}
\begin{tabular}{|m{2.8cm}|m{7cm}|m{3.8cm}|m{2cm}|}
\hline
& Sketch of the tube (the open parts are in grey) & main resonance mode &
frequencies \\
\hline 
\hline
flute, recorder, open organ pipe\ldots
 &
\epsfig{width=7cm,file=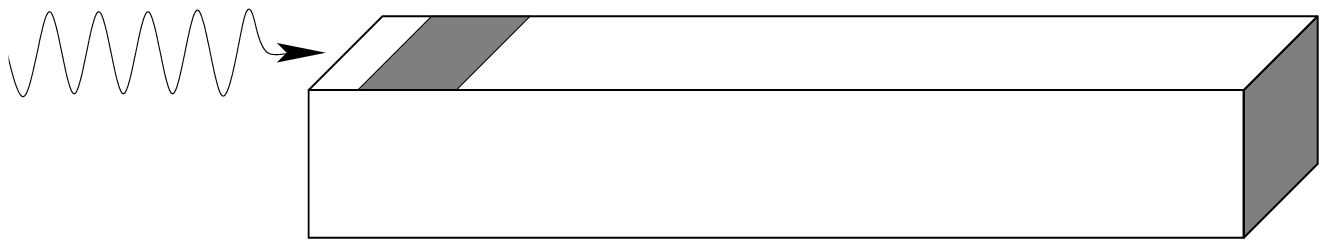}& \epsfig{width=4cm, height= 2cm,
file=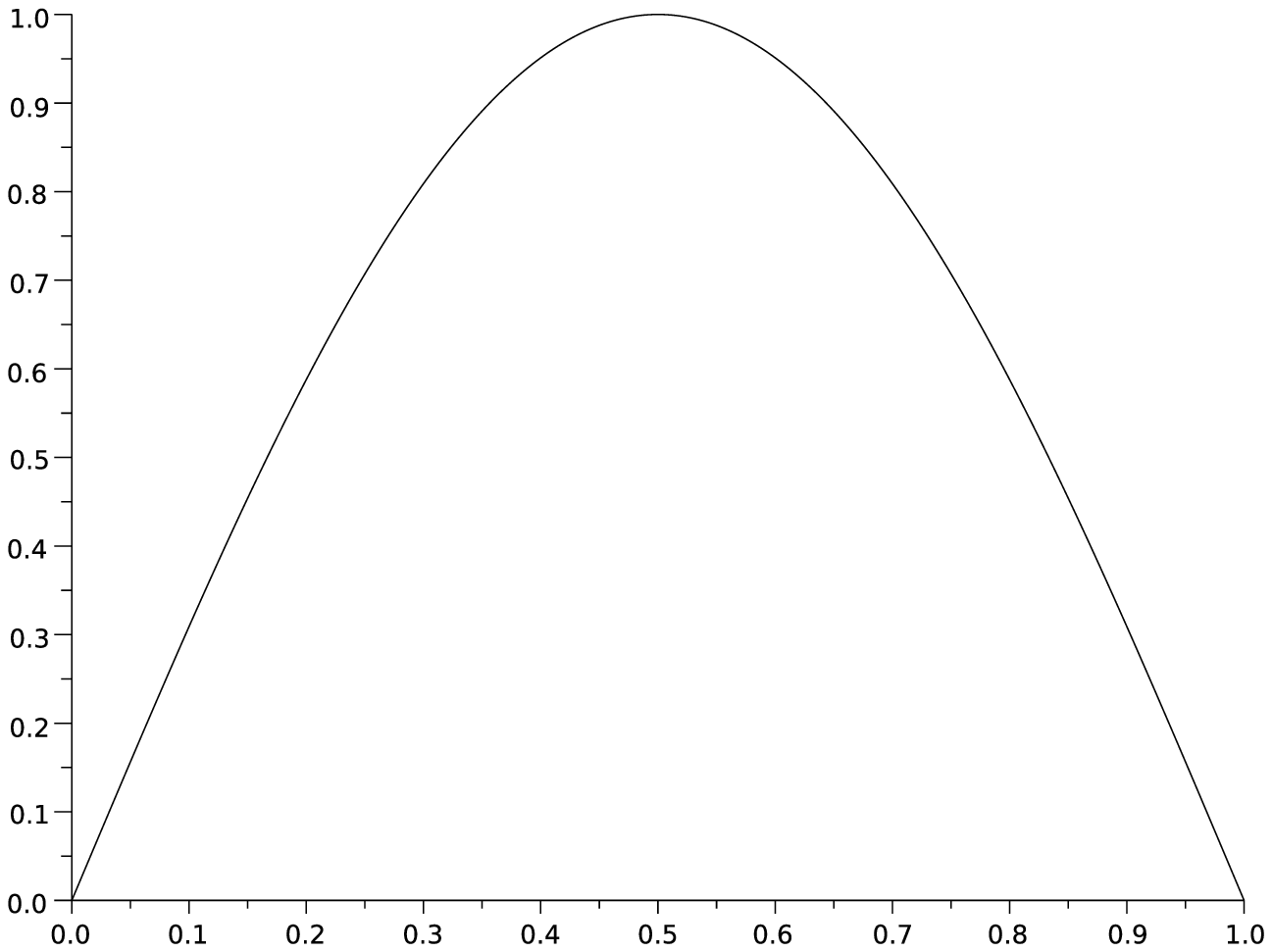}& $\pi/L$, $2\pi/L$, $3\pi/L$ \ldots\\
\hline
closed organ pipe, panpipes\ldots & \epsfig{width=7cm,file=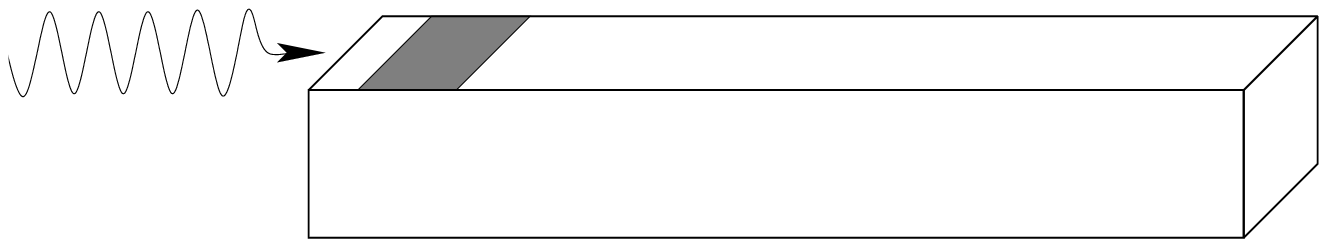}&
\epsfig{width=4cm,height= 2cm,file=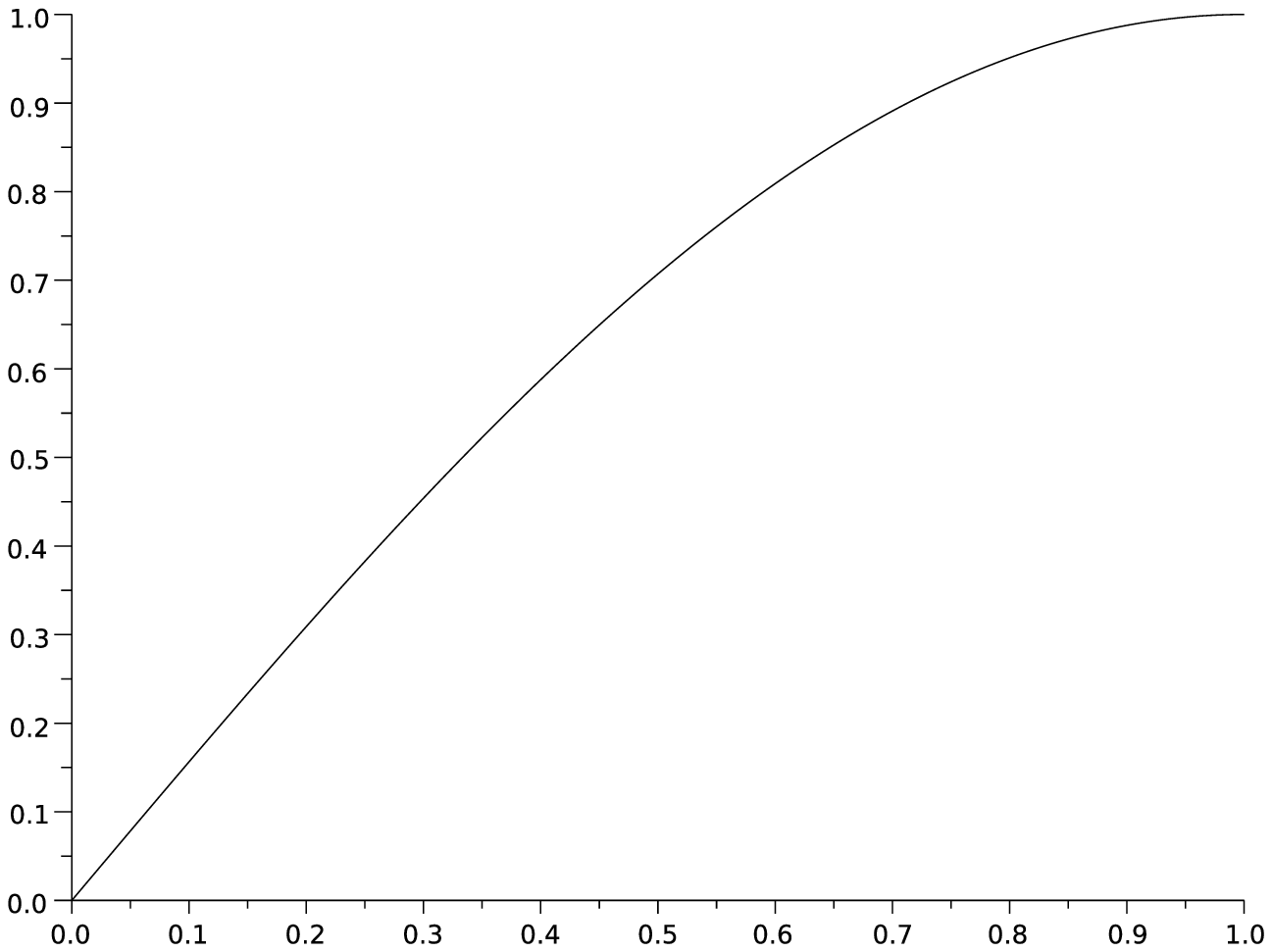}& $\pi/2L$, $3\pi/2L$,
$5\pi/2L$ \ldots\\
\hline
reed instruments (clarinet, oboe\ldots) & \epsfig{width=7cm,file=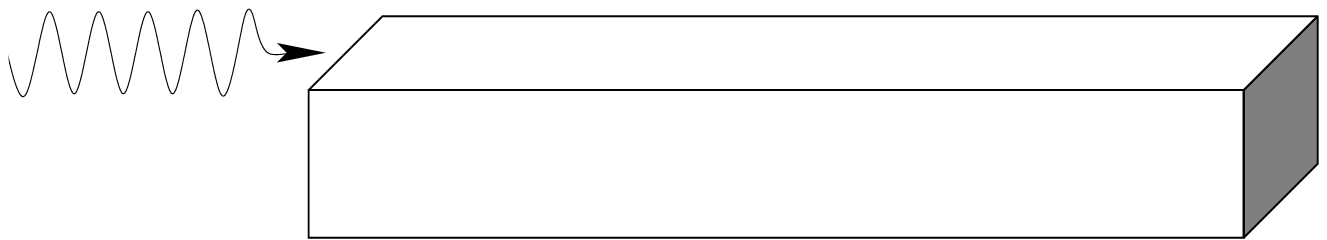}&
\epsfig{width=4cm,height= 2cm,file=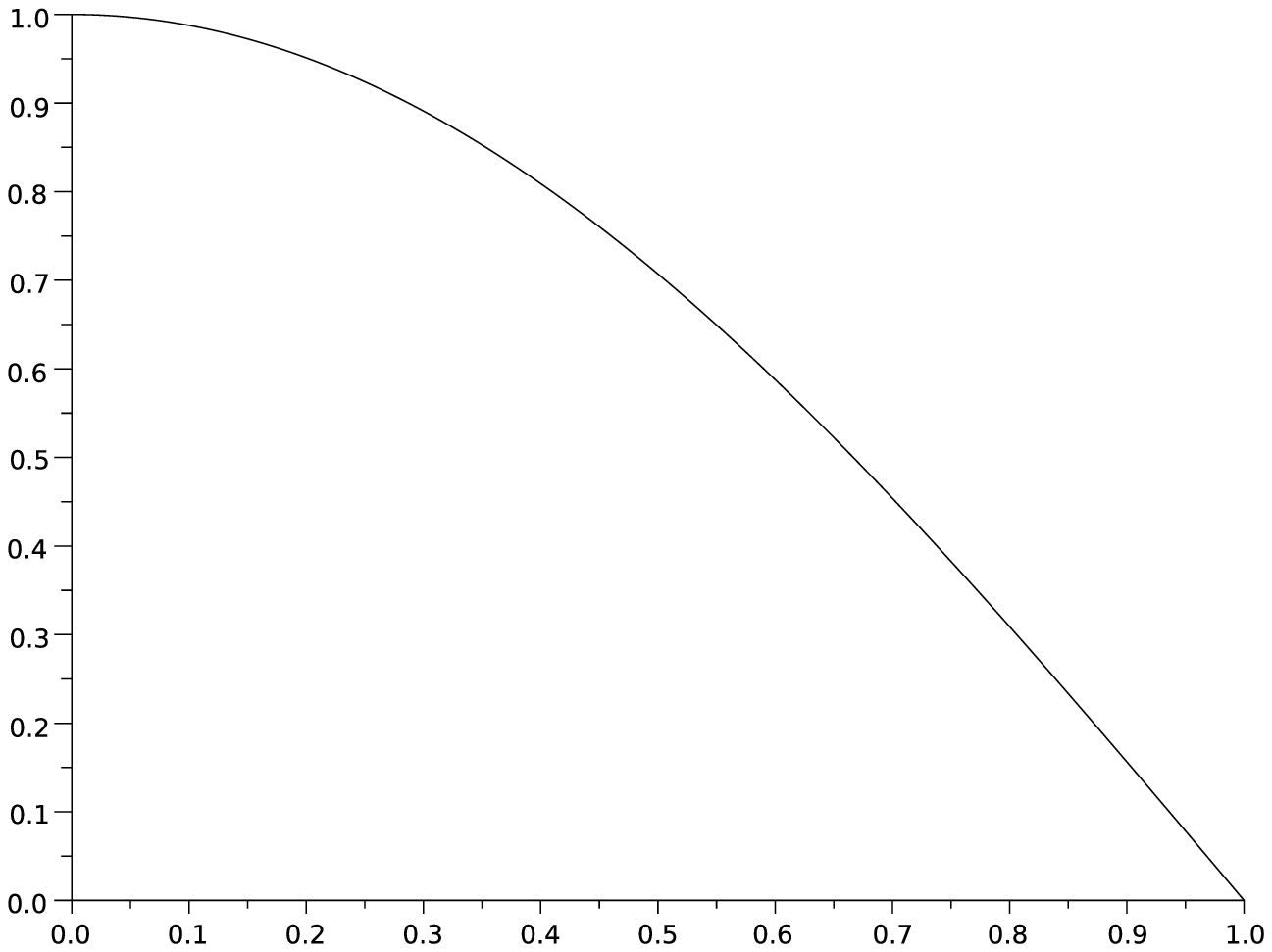}& $\pi/2L$, $3\pi/2L$,
$5\pi/2L$ \ldots\\
\hline
\end{tabular}
\end{center}
\caption{\it Three different kind of tubes without holes and their
resonances.} 
\label{tubes}
\end{figure}

We can roughly classify the tube of the wind instruments in three different
categories, depending on which end of the tube is open. See Figure \ref{tubes}.
It is known since a long time that the resonances of the tubes of
Figure \ref{tubes} can be approximated by the spectrum of the
one-dimensional Laplacian operator on $(0,L)$ with either Dirichlet or Neumann
boundary conditions, depending on whether  the corresponding end is open or not
(see for example \cite{CH}).
Notice that this rough approximation can explain simple facts: a tube with a
closed end sounds an octave lower than an open tube of the same length (enabling
for example to make shorter organ pipes for low notes) and moreover it produces
only the odd harmonics (explaining the particular sounds of reed instruments).

In this article, we study how the one-dimensional limit is affected by 
opening one of the holes of the flute, say a hole at position $a\in(0,L)$. At
first sight, one may think that it is equivalent to cutting the tube at the
place of the open hole. In other words, the note is the same as the one
produced by a tube of length $a$. This is
roughly true for flutes with large holes as the modern transverse flute, except
that one must add a small correction and the length $\tilde a$ of the
equivalent tube is slightly larger than $a$. This length $\tilde a$ is called
the effective length. This kind of approximation seems to be the most used one 
by acousticians. It states that the resonances of the tube with an open hole
are:
a fundamental frequency\footnote{We use in this article the mathematical habit
to identify
the frequencies to the eigenvalues of the wave operator. To obtain the real
frequencies corresponding to the sound of the flute, one has to divide them by
$2\pi$} $\pi/\tilde a$ and harmonics $k\pi/\tilde a$, $k\geq 2$.
However, the approximation of the resonances by the ones of a tube of length
$\tilde a$ is too rough for flutes with
small holes as the baroque flute or the recorder. In particular, the
approximation by effective length fails to explain the
following observations, for which we refer e.g. to \cite{Benade} and
\cite{Wolfe}:\\
- the effective length depends on the frequency of the waves in the tube. In
other words, the harmonics are not exact multiples of the fundamental
frequency.\\
- closing or opening one of the holes placed after the first open hole of
the tube changes the note of the flute. This enables to obtain some notes by
fork fingering, as it is common in baroque flute or recorder. We also enhance
that some effects of the baroque flute or of the recorder are produced by
half-holing, that is that by half opening a hole (some flutes have even holes
consisting in two small close holes to make half-holing easier). In these
cases, the effective length $\tilde a$ is not only related to the position $a$
of the first open hole, which makes the method of approximation by effective
length less relevant.\\

The purpose of this article is to obtain an explicit one-dimensional
mathematical model for the flute with a open hole, which could be more
relevant in the case of small holes than the approximation by effective
length. The models used by the acousticians are based on the notion of
impedance. The model introduced here rather uses the framework of
differential operators.

\vspace{3mm}

{\noindent \bf The thin domains techniques.}\\
The fact that the behaviour of thin three dimensional objects as a rope or
a plate can be approximated by one- or two-dimensional
equations has been known since a long time, see \cite{Hadamard} and \cite{CH}
for example. In general, a thin domain problem consists in a partial
differential equation $(E_\eps)$ defined in a domain $\Omega_\eps$ of dimension
$n$, which has $k$ dimensions of negligible size with respect to the other
$n-k$ dimensions. The aim is then to obtain an approximation of the problem by
an equation $(E)$ defined in a domain $\Omega$ of dimension $n-k$.
It seems that the first modern rigorous studies of such approximations mostly
date back to the late 80's: \cite{LS}, \cite{Anne}, \cite{AHH},
\cite{Hale-Raugel}, \cite{Schatzman}... There exists
an enormous quantity of papers dealing with thin domain problems of many
different types. We refer to \cite{Raugel} for a presentation of the subject and
some references.

In this paper, the domain $\Omega_\eps$ is the thin tube of the flute and we
hope to model the behaviour of the internal air pressure by a one-dimensional
equation. It is well known that the wave equation in a simple tube can be
approximated by the one-dimensional wave equation. Even the case of a far more
complicated domain squeezed along some dimension is well understood, see
\cite{PR} and the references therein. We will assume in this paper that the
open parts of the tube yield a Dirichlet boundary condition for the pressure in
the tube. In fact, we could study the whole system of a thin tube
connected to a large room and show that, at a first order of approximation, the
effect of the connection with a large domain is the same as a the one
of a Dirichlet
boundary condition, see \cite{Beale}, \cite{AHH}, \cite{JM} and the other works
related to the ``dumbbell shape'' model. The main difficulty of our problem
comes from the different scales: the open hole on the side of the
tube is of size $\eps^2$, whereas the diameter of the tube is of size $\eps$.
Thin domains involving different order of thickness have been studied in
\cite{Ciuperca}, \cite{Murat}, \cite{Murat-CRAS1}, \cite{Murat-CRAS2} and the
related works. The methods used in this paper are mainly based on these last
articles of J. Casado-D\'\i az, M. Luna-Laynez and F. Murat.

\vspace{3mm}

\begin{figure}[tp]
\begin{center}
\resizebox{0.9\textwidth}{!}{\input{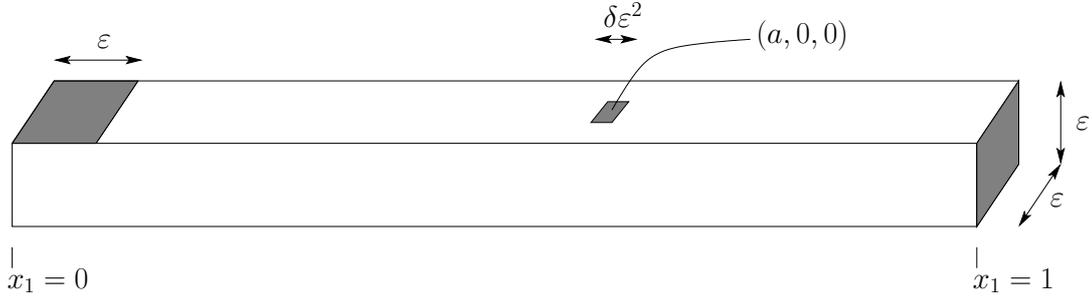}}
\end{center}
\caption{\it The domain $\Omega_\eps$. The grey parts correspond to Dirichlet
boundary conditions, the other ones to Neumann boundary conditions.} 
\label{domaine}
\end{figure}

{\noindent \bf Notations and main result.}\\
For $\eps>0$, we consider the domain
$$\Omega_\eps=(0,1)\times(-\eps,0)\times(-\eps/2,\eps/2)~.$$
We split any $x\in\Omega_\eps$ as $x=(x_1,x_2,x_3)=(x_1,\tilde x)$. Let $a\in
(0,1)$ and $\delta>0$. We denote by $\Delta_\eps$ the positive Laplacian
operator with the following boundary conditions:
$$\left\{\begin{array}{lll} \text{Dirichlet B.C.
}&~u=0~&\text{ on }~~~
(0,\eps)\times\{0\}\times(-\eps/2,\eps/2)\\
&& ~~~~~~~~ \cup \{1\}\times (-\eps,0)\times
(-\eps/2,\eps/2)\\ 
&& ~~~~~~~~ \cup (a-\delta\eps^2/2,a+\delta\eps^2/2) \times
\{0\} \times (-\delta\eps^2/2,\delta\eps^2/2)\\
\text{Neumann B.C. }&~\partial_\nu u=0~&\text{ elsewhere.}
\end{array}\right.$$
We denote by $H^1_0(\Omega_\eps)$ the Sobolev space corresponding to the above 
Dirichlet boundary conditions. The domain $\Omega_\eps$ is represented in Figure
\ref{domaine}.

In this paper, we show that, when $\eps$ goes to $0$, the spectrum of
the operator $\Delta_\eps$ converges to the one of the one-dimensional operator
$A$, defined by
$$A:\left(\begin{array}{ccc} D(A) & \longrightarrow & L^2(0,1)\\ u &
\longmapsto & -u'' \end{array} \right) $$
where $D(A)=\{u\in H^2((0,a)\cup(a,1)) \cap
H^1_0(0,1)~|~u'(a^+)-u'(a^-)=\alpha\delta u(a)\}$
and where $\alpha$ is the positive constant given by
\begin{equation}\label{def-alpha}
\alpha=\int_K |\grad \zeta|^2~,
\end{equation}
with $\zeta$ being the auxiliary function introduced in Proposition
\ref{prop-zeta} below.

Notice that both $\Delta_\eps$ and $A$ are positive definite self-adjoint
operators and that
\begin{equation}\label{eq-var-D}
\forall u,v\in D(\Delta_\eps)~,~~\langle
\Delta_\eps u|v\rangle_{L^2(\Omega_\eps)}=\int_{\Omega_\eps} \grad
u(x)\grad v(x)dx~,
\end{equation}
\begin{equation}\label{eq-var-A}
\forall u,v\in D(A)~,~~\langle
Au|v\rangle_{L^2(0,1)}=\int_0^1 u'(x) v'(x)dx+\alpha \delta u(a)v(a)~.
\end{equation}
Let $0<\lambda^1_\eps<\lambda^2_\eps \leq \lambda^3_\eps\leq\ldots$ be the
eigenvalues of $\Delta_\eps$ and let $0<\lambda^1<\lambda^2
\leq \lambda^3\leq\ldots$ be the ones of $A$. The purpose of this paper is
to prove the following result.
\begin{theorem}\label{th}
When $\eps$ goes to $0$, the spectrum of $\Delta_\eps$ converges to the one of
$A$ in the sense that
$$ \forall
k\in\NN^*~,~~\lambda^k_\eps\xrightarrow[~~~\eps\longrightarrow 0~~~]{}
\lambda^k~.$$
\end{theorem}

Theorem \ref{th} yields a new model for the flute, which is discussed in
Section \ref{discuss}. The proof of Theorem \ref{th} consists in showing lower-
and upper-semicontinuity of the spectrum, which is done is Sections \ref{low}
and \ref{up} respectively. We use scaling techniques consisting in focusing to
the hole at the place $(a,0,0)$. These techniques follow the ideas of
\cite{Murat} (see also \cite{Murat-CRAS1} and \cite{Murat-CRAS2}). The
corresponding technical background is introduced in Section \ref{zoom}.

\vspace{3mm}

\noindent {\bf Acknowledgements:} the interest of the author for the
mathematical models of flutes started with a question of Brigitte Bid\'egaray
and he discovered the work of J. Casado-D\'\i az, M. Luna-Laynez and F. Murat
following a discussion with Eric Dumas. The author also thanks the referee for
having reviewed this paper so carefully and so quickly.


\section{Discussion}\label{discuss}

\begin{figure}[tp]
\begin{center}
{\epsfig{width=6.5cm,file=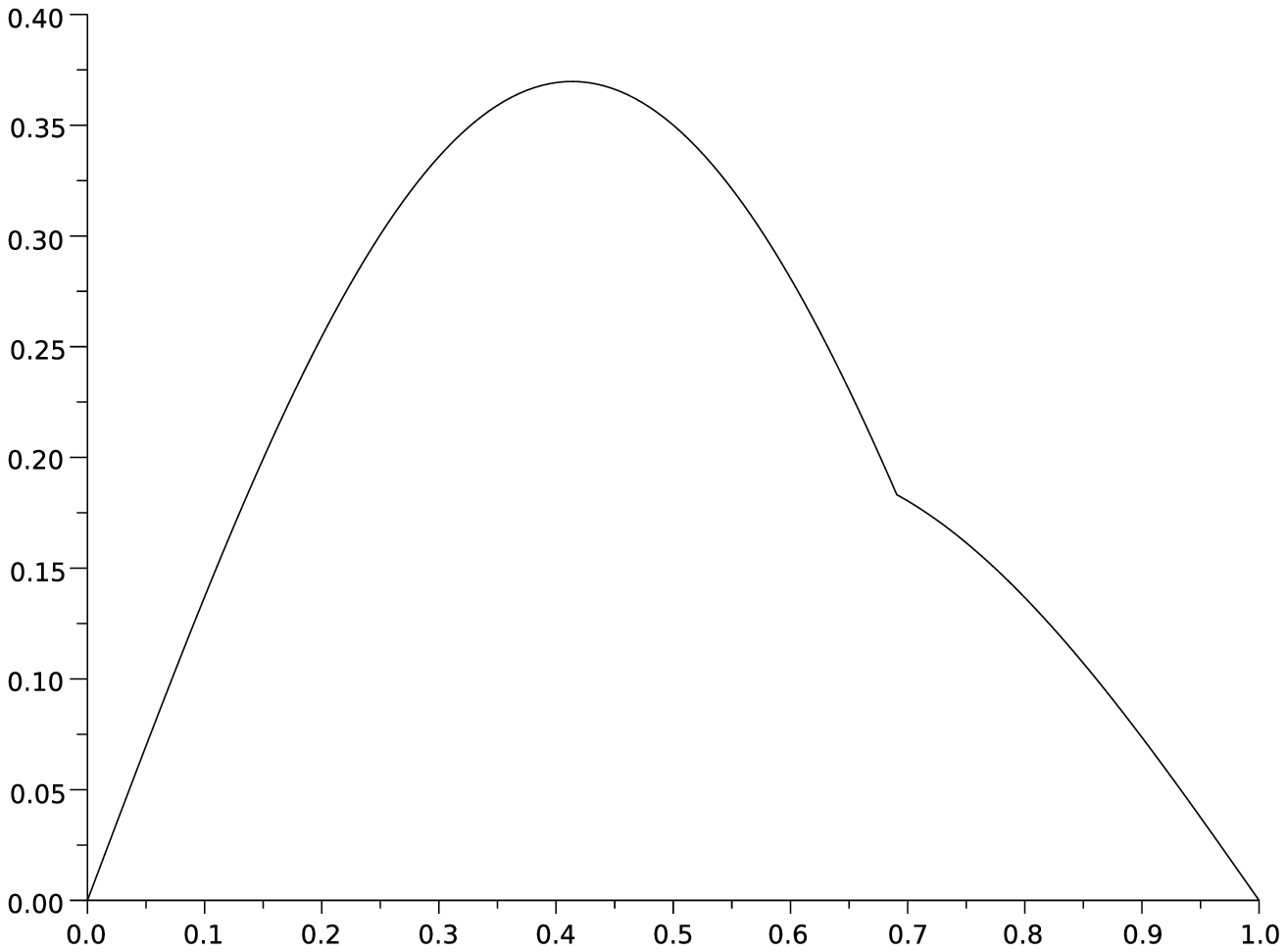}}
{\epsfig{width=9.5cm,height=5cm,file=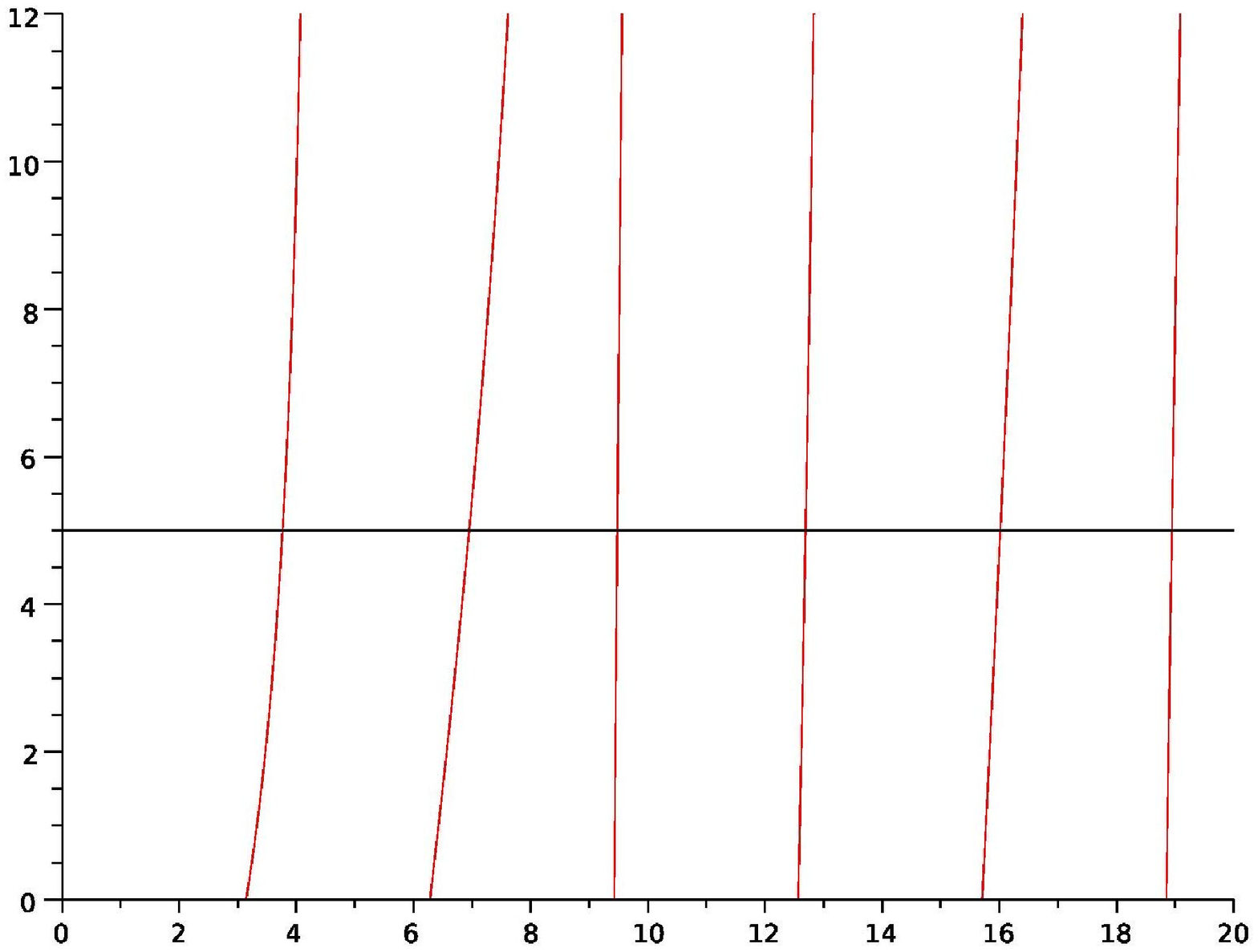}}
\end{center}
\caption{\it Left, the first eigenfunction of $A$, i.e. the fundamental mode of
resonance of the flute with an open hole. Right, the graphic of the function
$f_a$ and the intersections with the line $y=\alpha\delta$ giving the
frequencies of the flute. The references values are $a=0.7$ and
$\alpha\delta=5$.}
\label{fct}
\end{figure}

First, let us compute the frequencies of the flute with an open hole,
following the model yielded by Theorem \ref{th}.  
Theorem \ref{th} deals with the spectrum of $\Delta_\eps$, whereas the
resonances of the pressure in a flute follow the wave equation
$\partial^2_{tt}u=-\Delta_\eps u$ (remind that $\Delta_\eps$ denotes the
positive Laplacian operator). Therefore, the relevant eigenvalues are in fact
the ones of the operator {\footnotesize $\left(\begin{array}{cc}0&Id\\
-\Delta_\eps& 0 \end{array}\right)$} which are $\pm i \sqrt{\lambda_\eps^k}$.
Theorem \ref{th} shows that the frequencies $\sqrt{\lambda_\eps^k}$ are
asymptotically equal to the frequencies $\mu>0$ such that
$\mu^2$ is an eigenvalue of $A$. A straightforward computation shows that
$\mu^2$ is an eigenvalue of $A$, with corresponding eigenfunction $u$, if and
only if 
$$u(x)=\left\{\begin{array}{ll} C \sin(\mu x)&~~x\in(0,a)\\
            C\frac{\sin(\mu a)}{\sin(\mu(1-a))} \sin(\mu(1-x))&~~x\in (a,1)
\end{array}\right. $$
with some $C\neq 0$ and with $\mu>0$ solving 
\begin{equation}\label{eq-mu}
\alpha\delta=\frac{-\mu\sin\mu}{\sin(\mu a)\sin(\mu(1-a))}:=f_a(\mu)~,
\end{equation}
see Figure \ref{fct}.

\vspace{3mm}
 
\noindent Using the above computations, we can do several remarks about the
resonances of
the flute with a small open hole, as predicted by our model.
\begin{itemize}
\item The eigenfunctions of $A$ corresponds to the expected profile of
the pressure in the flute with an open hole, see Figure \ref{fct} and the ones
of \cite{Bolton}, \cite{Coltman} and \cite{Wolfe}.
\item The note of the flute corresponds to the fundamental frequency
$\mu=\sqrt{\lambda^1}$. To obtain a given note, one can adjust both $\delta$
(the size of the hole) and $a$ (the place of the hole). This enables to place
smartly the different holes to obtain some notes by combining the opening of
several holes (fork fingering). We can also compute the change of frequency
produced by only half opening the hole (half-holing). Notice that changing the
shape of the hole affects the coefficient $\alpha$.
\item The overtones of the flute correspond to the other frequencies
$\mu=\sqrt{\lambda^k}$ with $k\geq 2$. We can see in Figure \ref{fct} that they
are not exactly harmonic, i.e. they are not multiples of the fundamental
frequency. This explains why the sound of flutes, which have only a small hole
opened, is uneven and not as pure as the sound produced by a simple tube. In
other words, our model directly explains the observation that the effective
length approximation depends on the considered frequency. Moreover, 
when $\mu$ increases, the slope of $f_a$ becomes steeper due to the factor
$\mu$ in \eqref{eq-mu} and the solutions of
\eqref{eq-mu} are closer to $\mu=k\pi$. This is consistant with the observation
that high frequencies are less affected by the presence of the hole than low
frequencies, see \cite{Wolfe} or \cite{Wolfe-web}. However, notice that this is
only roughly true since for example one can see on Figure \ref{fct} that the
second overtone is almost equal to $3\pi$, whereas the fourth one is less close
to $5\pi$. This comes from the fact that $a=0.7$ is almost a node of the mode
$\sin(3\pi x)$.
\item Of course, when $\delta=0$, we recover the equation $\sin\mu=0$
corresponding to the eigenvalue of the open tube without hole. When
$\delta\rightarrow +\infty$, i.e. when the hole is very large, we recover the
equations $\sin(\mu a)=0$ or $\sin(\mu(1-a))=0$, which correspond to two
separated tubes of lengths $a$ and $1-a$ (in fact the part $(a,1)$ is not
important because this is not the part of the tube which is excited by the
fipple). When the hole is of intermediate size, the fundamental frequency
corresponds to a tube of intermediate length $\tilde a\in(a,1)$, but the
overtones are not the same as the ones of the tube of length $\tilde a$.
\item The thin domain techniques used here are general and do not depend on the
fact that the section of the tube $\Omega_\eps$ is a square and not a disk. If
the surface $g(x)$ of the section of the tube is not constant (think at the
end of a clarinet), then the operator $\partial_{xx}^2$ in the definition of $A$
must be replaced by ${\small \frac1{g(x)}}\partial_x(g(x)\partial_x.)$, see
\cite{Hale-Raugel}. Of course, if there are several open holes, then other terms
of the type $\alpha \delta u(a)v(a)$ appear in \eqref{eq-var-A}.
\end{itemize}
To conclude, we obtain in this article a mathematical model for the flute with a
small open hole, which consists in a one-dimensional operator different from a
simple Laplacian operator. It yields simple explanation of some observations as
the fact that the overtones are not harmonic.


\section{Focusing on the hole: the rescaled problem}\label{zoom}

\begin{figure}
 \begin{center}
 {\resizebox{0.95\textwidth}{!}{\input{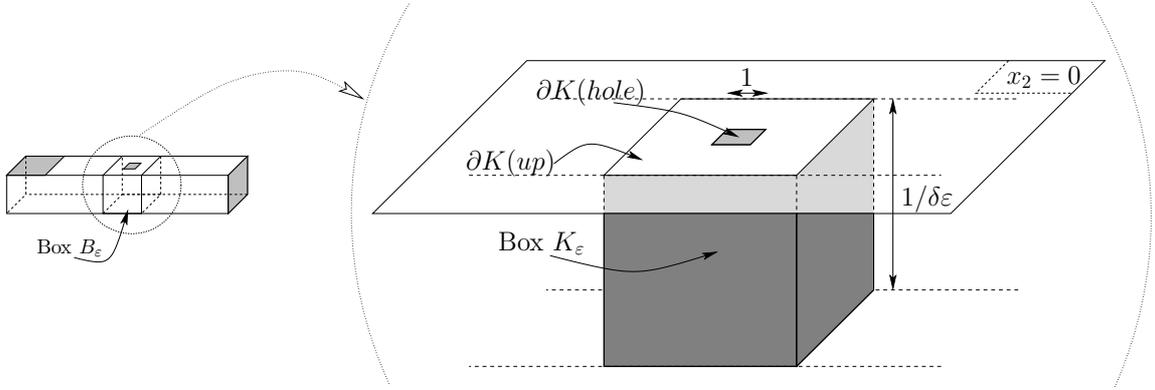}}}
 \end{center}
 \caption{\it The cube $K_\eps$, part of the half-space
$K=\{x\in\RR^3,~x_2<0\}$, and the corresponding boundaries. When $\eps$ goes to
$0$, the cube $K_\eps$ converges to $K$, whereas the hole $\partial K(hole)$
remains unchanged.}
 \label{K}
 \end{figure}

When $\eps$ goes to zero, if one rescales the domain $\Omega_\eps$ with a
ratio $1/(\delta\eps^2)$ to focus on the hole, then one sees the rescaling
domain $\Omega_\eps$ converging to the half-space $x_2<0$ (see Figure~\ref{K}).
The purpose of this section is to introduce the technical background to be able
study our problem in this rescaled frame. For the reader interested in
more details about the Poisson problem in unbounded domain, we refer to
\cite{Simader-Sohr}.

\subsection{The space $\dot H^1(K)$}
Let $K$ be the half-space $\{x\in\RR^3,~x_2<0\}$. For any $\eps>0$, we introduce
the cube 
$$K_\eps=\left(\frac{-1}{2\delta\eps},\frac{1}{2\delta\eps}
\right)\times\left(\frac{-1}{\delta\eps},0\right)\times
\left(\frac{-1}{2\delta\eps},\frac{1}{2\delta\eps}
\right)$$ 
as shown in Figure \ref{K}.
We denote by $\partial K_\eps(hole)$ the part of the boundary
$(-1/2,1/2)\times\{0\}\times (-1/2,1/2)$ corresponding to the
hole. We denote $\partial K_\eps(up)$ the remaining part of the upper face.
We also denote by $\partial K(hole)$ and $\partial
K(up)$ the corresponding parts of the boundary of the half-space $K$. See
Figure \ref{K}.

We introduce the space $\dot H^1(K)$ defined by
\begin{equation}\label{def-H-point}
\dot H^1(K)=\{ v \in H^1_{loc}(K)~,~\grad v\in L^2(K)\text{ and }v=0\text{ on
}\partial K(hole) \}~
\end{equation}
and we equip it with the scalar product 
\begin{equation}\label{def-scalar}
\langle \varphi|\psi\rangle_{\dot H^1(K)} =
\int_K\grad\varphi.\grad\psi~.
\end{equation}
We also introduce the space $\dot H^1_0(K)$ which is the completion of
\begin{equation}\label{def-H-point-0}
\Cc^\infty_0(\overline K)=\{\varphi \in\Cc^\infty(\overline
K)~/~\supp(\varphi)\text{ is compact
and }\varphi\equiv 0\text{ on }\partial K(hole)\} 
\end{equation}
with respect to the $\dot H^1$ scalar product defined in \eqref{def-scalar}.

Let $\chi\in\Cc^\infty(\overline K)$ be such that $\chi\equiv
1$ outside a compact set, $\chi\equiv 0$ on $\partial K(hole)$ and
$\partial_\nu\chi\equiv 0$ on $\partial K(up)$. Following \cite{Simader-Sohr},
we get the following results.
\begin{theorem}\label{th-spaces}
The spaces $\dot H^1(K)$ and $\dot H^1_0(K)$ equipped with the scalar product
\eqref{def-scalar} are Hilbert spaces and  
\begin{equation}\label{eq-th-spaces}
\dot H^1(K)=\dot H^1_0(K) \oplus  \RR \chi~,
\end{equation}
this sum being a direct sum of closed subspaces.

Moreover, a function $u\in \dot H^1(K)$ belongs to $\dot H^1_0(K)$ if and only
if it belongs to $L^6(K)$. As a consequence, the splitting of $u\in \dot
H^1(K)$ given by \eqref{eq-th-spaces} is uniquely determined by $u=\dot u +
\overline u \chi$, where 
$$\overline u = \lim_{\eps\rightarrow 0}  \frac 1{|K_\eps|}
\int_{K_\eps} u(x)dx$$
is the average of $u$, which is well defined.
\end{theorem}
\begin{demo}
The direct sum \eqref{eq-th-spaces} is a particular case of Theorem 2.15
of \cite{Simader-Sohr}. The equivalence between $u\in \dot H^1_0(K)$ and $u\in
L^6(K)$ is given by Theorem 2.8 of \cite{Simader-Sohr}. Let $u=\dot u +
c \chi$ with $\dot u \in \dot H^1_0(K)$ and $c\in\RR$. Since $\dot u\in
L^6(K)$, we have $\int_{K_\eps} |\dot u| \leq C |K_\eps|^{5/6}$ and thus the
average of $\dot u$ is well defined and equal to $0$. Since the average of
$\chi$ is well defined and equal to $1$, the average of $u$ is also well
defined and it is equal to $c$.
\end{demo}

\subsection{The function $\zeta$}
We now introduce the function $\zeta$, which is used to define the
coefficient $\alpha$ in \eqref{def-alpha}.
\begin{prop}\label{prop-zeta}
There is a unique weak solution $\zeta$ of 
\begin{equation}\label{eq-zeta}
 \left\{\begin{array}{ll}
         \Delta \zeta =0&\text{ on }K\\
         \zeta=0&\text{ on }\partial K(hole)\\
         \partial_\nu \zeta=0&\text{ on }\partial K(up)\\
        \overline\zeta=1&
        \end{array}
\right.
\end{equation}
in the sense that $\zeta \in \dot H^1(K)$, $\overline \zeta=1$ and 
$$\forall \varphi\in\Cc^\infty_0(K)~,~~\int_K \grad \zeta.\grad\varphi=0~. $$
\end{prop}
\begin{demo}
Theorem \ref{th-spaces} shows that $\zeta=\chi+\dot\zeta$ with $\dot\zeta\in
\dot H^1_0(K)$. Then, Proposition \ref{prop-zeta} is a direct application of
Lax-Milgram Theorem to the variational equation
$$\forall \dot\varphi\in  \dot H^1_0(K)~,~~\int_K \grad
\dot\zeta\grad\dot\varphi = - \int_K \Delta\chi.\dot\varphi~.$$
See \cite{Simader-Sohr}
for a discussion on this kind of variational problems.
\end{demo}

The function $\zeta$ yields a different way to write the scalar product in $\dot
H^1(K)$.
\begin{prop}\label{prop-redef-scalar}
The function $\zeta$ is the orthogonal projection of $\chi$
on the orthogonal space of $\dot H^1_0(K)$ in $\dot H^1(K)$. 

Thus, for all $u$ and $v$ in $\dot H^1(K)$, there exist two unique functions
$\dot u$ and $\dot v$ in $\dot H^1_0(K)$ such that $u=\dot u +\overline u \zeta$
and $v=\dot v +\overline v \zeta$. Moreover, 
$$\langle u | v \rangle_{\dot H^1(K)} ~=~ \int_K \grad \dot u (x). \grad\dot
v(x)~dx~+~\alpha \overline u .\overline v~,$$
where $\alpha$ is defined by \eqref{def-alpha}.
\end{prop}

\subsection{Weak $K_\eps-$convergence}
As one can see in Figure \ref{K}, if $(u_\eps)$ is a sequence of functions
defined in $\Omega_\eps$, then the rescaled functions $w_\eps$ are only defined
in the box $K_\eps$ and not in the whole space $K$. Hence, we have to introduce
a suitable notion of weak convergence. 
\begin{prop}\label{prop-weak}
Let $(w_\eps)_{\eps>0}$ be a sequence of functions of $H^1(K_\eps)$ vanishing
on $\partial K_\eps(hole)$. Assume that 
$$\exists C>0~,~~\forall \eps>0~,~~\int_{K_\eps} |\grad w_\eps|^2 \leq C~.$$
Then, there exists a subsequence $(w_{\eps_n})_{n\in\NN}$, with
$\eps_n\rightarrow 0$, which converges weakly to a function $w_0\in\dot H^1(K)$
in the sense that 
$$\forall \varphi\in \dot H^1(K)~,~~\int_{K_{\eps_n}} \grad
w_{\eps_n}\grad\varphi \xrightarrow[~~\eps_n\longrightarrow
0~~]{} \int_K \grad w_0 \grad \varphi~.$$
Moreover, the average of $w_0$ is given by 
\begin{equation}\label{eq-prop-weak}
\overline w_0=\lim_{ \eps_n\rightarrow 0} \frac 1{|K_{\eps_n}|}
\int_{K_{\eps_n}} w_\eps~.
\end{equation}
\end{prop}
Before starting to prove Proposition \ref{prop-weak}, we recall 
Poincar\'e-Wirtinger inequality. 
\begin{lemma}\label{prop-Poincare}
(Poincar\'e-Wirtinger inequality)\\
There exists a constant $C>0$ such that, for any $\eps>0$ and any function
$\varphi\in H^1(K_\eps)$, 
\begin{equation}\label{eq-poincare}
\int_{K_\eps} \left|\varphi(x)- \frac 1{|K_\eps|}
\int_{K_\eps} \varphi(s)ds\right|^6dx \leq C \left(\int_{K_\eps}
|\grad \varphi(x)|^2dx\right)^3~.
\end{equation}
\end{lemma}
\begin{demo}
First, let us set $\eps=1$. The classical Poincar\'e inequality (see
\cite{Evans} for example) states that 
$$\int_{K_1} \left|\varphi(x)- \frac 1{|K_1|}
\int_{K_1} \varphi(s)ds\right|^2dx \leq C \int_{K_1}
|\grad \varphi(x)|^2dx~.$$
Thus, the right-hand side controls the $H^1-$norm of $\varphi-\overline\varphi$.
Then, the Sobolev inequalities shows that \eqref{eq-poincare} holds for
$\eps=1$. Now, the crucial point is to notice that the constant
$C$ in \eqref{eq-poincare} is independent of the size of the cube $K_\eps$ since
both sides of the inequality behave similarly with respect to scaling.
\end{demo}

{ \noindent \emph{\textbf{Proof of Proposition \ref{prop-weak}~:}}} 
First notice that $\dot H^1_0(K)$ is separable due to the density of
$\Cc^\infty_0-$functions. Hence, $\dot H^1(K)$ is also separable and by a
diagonal extraction argument, we can extract a subsequence $\eps_n\rightarrow 0$
such that for all $\varphi\in \dot H^1(K)$, $\int_{K_{\eps_n}} \grad
w_{\eps_n}\grad\varphi$ converges to a limit $L(\varphi)$ with $L(\varphi)\leq
C\|\varphi\|_{\dot H^1}$. By Riesz representation Theorem, there exists
$w_0\in \dot H^1(K)$ such that $L(\varphi)=\langle w_0 |\varphi \rangle$.

To prove \eqref{eq-prop-weak}, we follow the arguments of \cite{Murat}. We set
$\overline w_\eps=\frac 1{|K_{\eps}|}\int_{K_{\eps}} w_\eps$. Let
$p\in\NN$. Lemma \ref{prop-Poincare} and the fact that $\int_{K_\eps} |\grad
w_\eps|^2$ is bounded, show that there exists a constant $C$ independent of
$\eps$ such that $\int_{K_\eps} |w_\eps(x)- \overline w_\eps|^6dx \leq C$.
Thus,  
\begin{equation}\label{eq-prop-weak-3}
\forall \eps< \frac 1p ~,~~ \int_{K_{1/p}} |w_\eps(x)- \overline w_\eps|^6dx
\leq C~.
\end{equation}
By Sobolev inequality, we know that $w_\eps$ is bounded in $L^6(K_{1/p})$
(remember that $w_\eps$ vanishes on $\partial K_{1/p}(hole)$). Thus
$\overline w_\eps$ is bounded and up to extracting another subsequence, we can
assume that $\overline w_{\eps_n}$ converges to some limit $\beta\in\RR$. By a
diagonal extraction argument, we can
also assume that $w_{\eps_n}$ converges to $w_0$ weakly in $L^6(K_{1/p})$, for
any $p\in\NN$. As a consequence, \eqref{eq-prop-weak-3} implies that 
$$\int_{K_{1/p}} |w_0(x)- \beta|^6dx \leq \limsup_{\eps\longrightarrow 0}
\int_{K_{1/p}} |w_\eps(x)- \overline w_\eps|^6dx \leq C~.$$
Since the previous estimate is uniform with respect to $p\in\NN$ and since
$K_{1/p}$ grows to $K$ when $p$ goes to $+\infty$, we obtain that
$w_0-\beta$ belongs to $L^6(K)$ and thus $w_0-\beta\chi \in L^6(K)$. Theorem
\ref{th-spaces} shows that $w_0-\beta\chi$ belongs to $\dot H^1_0(K)$ i.e.
$\beta=\overline w_0$. {\hfill$\square$\\   \vspace{0.4cm}}


\section{Lower-semicontinuity of the spectrum}\label{low}

This section is devoted to the following result.
\begin{prop}\label{prop1}
$$\forall k\in\NN^*~,~~0\leq\displaystyle \limsup_{\eps\rightarrow 0}
\lambda^k_\eps\leq \lambda^k~.$$
\end{prop}
\begin{demo}
Let $(u^k)$ be a sequence of eigenfunctions of $A$ corresponding to
the eigenvalues $\lambda^k$. Since $A$ is symmetric, we may assume that
$\langle u^k|u^j\rangle_{L^2(0,1)}=0$ for $k\neq j$. The main idea of the proof
of Proposition \ref{prop1} is to construct an embedding 
$I_\eps:H^1_0(0,1)\rightarrow H^1_0(\Omega_\eps)$ such that the functions
$I_\eps u^k$ are almost $L^2-$orthogonal and such that
\begin{equation}\label{eq-demo-prop1}
\frac{\int_{\Omega_\eps} |\grad I_\eps u^k|^2}{\int_{\Omega_\eps} |
I_\eps u^k|^2} \xrightarrow[~~~\eps\longrightarrow 0 ~~~]{} \lambda^k~.
\end{equation}
The definition of the embedding 
$I_\eps:H^1_0(0,1)\rightarrow H^1_0(\Omega_\eps)$ is as follows.

\noindent \underline{Far from the hole:} we split the functions
$u^k$ into two parts $u^k_{|(0,a)}$ and $u^k_{|(a,1)}$, we slightly rescale them
so that they are defined in $(\eps,a-\eps/2)$ and $(a+\eps/2,1)$
respectively, and we embed both parts in $L^2(\Omega_\eps)$ by setting
$$\varphi_\eps^k(x)=u^k\left(\frac
a{a-3\eps/2}(x_1-\eps)\right)~~\text{ and
}~~\psi_\eps^k(x)=u^k\left(a+\frac{1-a}{
1-a-\eps/2}(x_1-a-\eps/2)\right)~.$$

\noindent \underline{Near the hole:} let $\zeta\in \dot H^1(K)$ be as in
Proposition \ref{prop-zeta} and let $\zeta=\dot \zeta+\chi$ be the splitting
given by Theorem \ref{th-spaces} (where we use that $\overline \zeta=1$ by
definition). By the definition of $\dot H^1_0(K)$, there exists a
sequence of functions $(\dot\zeta_\eps)\in \Cc^\infty_0(K)$ converging to $\dot
\zeta$ in $\dot H^1_0(K)$. Therefore, there exists a sequence
$\zeta_\eps=\dot\zeta_\eps+\chi \in \Cc^\infty(K)\cap \dot H^1(K)$
such that $\zeta_\eps\equiv 1$ outside a compact set and
$(\zeta_\eps)$ converges strongly to $\zeta$ when $\eps$ goes to zero. Notice
that we may assume that $\zeta_\eps\equiv 1$ outside a compact set of the cube
$K_\eps$ defined in Section \ref{zoom}. We set
$\tilde\zeta_\eps(x)=\zeta_{\eps}((x-(a,0,0))/(\delta\eps^2))$ and $I_\eps
u^k=u^k(a)\tilde\zeta_\eps$ in the cube $B_\eps=(a,0,0)+\delta\eps^2 K_\eps$.

\noindent \underline{Summarizin}g\underline{:} the whole
embedding $I_\eps$ is described by Figure
\ref{Ieps}.
\begin{figure}[tp]
\begin{center}
\resizebox{0.8\textwidth}{!}{\input{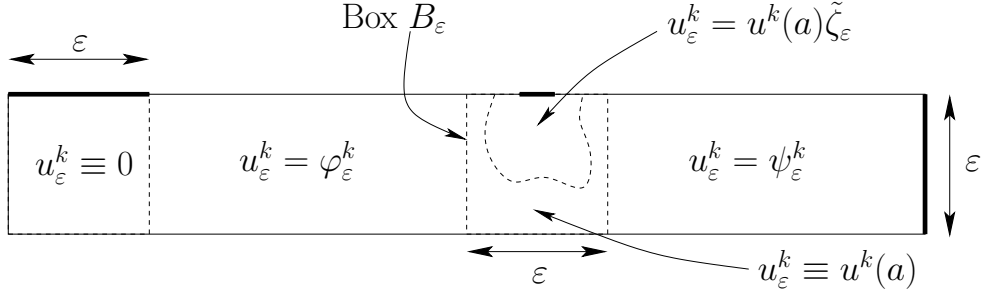}}
\end{center}
\caption{\it The embedding $u^k_\eps=I_\eps u^k$ of $u^k\in
H^1_0(0,1)$ into $H^1_0(\Omega_\eps)$ (lateral view).}
\label{Ieps}
\end{figure}

\noindent \underline{Calculin}g\underline{ the scalar
}p\underline{roducts:} 
by change of variables, we have 
$$\int_{B_\eps} |\tilde\zeta_{\eps}|^2=\delta^3
\eps^6 \int_{K_\eps} |\zeta_{\eps}|^2\leq \delta^3
\eps^6 \left(\left( \int_{K_\eps} |\dot \zeta_{\eps}|^2 \right)^{1/2}+ \left(
\int_{K_\eps} |\chi|^2 \right)^{1/2} \right)^2~.$$ 
Since $\chi$ is a bounded $\Cc^\infty$ function and since the volume of $K_\eps$
is of order $1/\eps^3$, we have $\int_{K_\eps} |\chi|^2=O(1/\eps^3)$.
Due to Theorem \ref{th-spaces}, $\dot \zeta_{\eps}$ converges to $\dot\zeta$ in
$L^6(K)$ and thus
$$\int_{K_\eps} |\dot \zeta_{\eps}|^2\leq \left(\int_{K_\eps} |\dot
\zeta_{\eps}|^6\right)^{1/3} \left(\int_{K_\eps} 1 \right)^{2/3}=O(1/\eps^2)~.
$$
Therefore, we get that $\int_{B_\eps} |\tilde\zeta_{\eps}|^2=O(\eps^3)$.
Thus, the $L^2-$norm of $u^k_\eps=I_\eps u^k$ is mostly due to the $L^2-$norms
of $\varphi_\eps^k$ and $\psi_\eps^k$ and so, for any $k$ and $j$,
\begin{equation}\label{eq-demo-prop1-2}
\langle u_\eps^k |u_\eps^j\rangle_{L^2(\Omega_\eps)}=\eps^2 \langle u^k
|u^j\rangle_{L^2(0,1)}+o(\eps^2)~.
\end{equation}
On the other hand, we have $\int_{B_\eps} |\grad \tilde\zeta_{\eps}|^2=\delta
\eps^2 \int_{K_\eps} |\grad \zeta_{\eps}|^2$. Since $(\zeta_\eps)$ converges to
$\zeta$ in $\dot H^1(K)$ and due to the definition \eqref{def-alpha} of
$\alpha$, $\int_{K_\eps} |\grad \zeta_{\eps}|^2$ converges to $\alpha$.
Therefore, for any $k$ and $j$,
\begin{align}
\int_{\Omega_\eps} \grad  u_\eps^k \grad  u_\eps^j  &=\int_{\Omega_\eps} \grad 
\varphi_\eps^k\grad 
\varphi_\eps^j+ \int_{\Omega_\eps} \grad 
\psi_\eps^k\grad \psi_\eps^j+u^k(a)u^j(a)\int_{B_\eps} |\grad \tilde
\zeta_{\eps}|^2\nonumber\\
&= \eps^2  \left(\int_0^a
\partial_x u^k(x)\partial_x u^j(x)  dx+\int_a^1
\partial_x u^k(x) \partial_x u^j(x) dx\right.\nonumber\\  
&~~~~~~~~~~~~~+\left.\delta u^k(a)u^j(a)
\int_{K_\eps} |\grad \zeta_{\eps}|^2 \right) + o(\eps^2)~\nonumber\\
&= \eps^2  \left(\int_0^1
\partial_x u^k(x) \partial_x u^j(x)dx+\alpha\delta u^k(a) u^j(a)  \right) +
o(\eps^2)~\nonumber \\
&= \eps^2 \langle A u^k|u^j\rangle_{L^2} + o(\eps^2)~.\label{eq-demo-prop1-3}
\end{align}
Hence the previous estimates yield the limit \eqref{eq-demo-prop1}.

\noindent \underline{A}pp\underline{l}y\underline{in}g\underline{
the Min-Max formula:} for $\eps$ small enough, \eqref{eq-demo-prop1-2} implies
that the functions $u_\eps^k$ are linearly independent. Due to the Min-Max
Principle (see \cite{Reed-Simon} for example), we know that
\begin{equation}\label{eq-prop1-1}
\lambda^k_\eps ~\leq ~ \min_{p_1<p_2<\ldots<p_k}~ \max_{c\in\RR^k_*}~
\frac{\int_{\Omega_\eps} |\sum_{i=1}^k c_i\grad I_\eps
u^{p_i}|^2}{\int_{\Omega_\eps} |\sum_{i=1}^k c_i I_\eps u^{p_i}|^2}~.
\end{equation}
The above estimates \eqref{eq-demo-prop1-2} and \eqref{eq-demo-prop1-3} show
that, for any $c\in\RR^k_*$, we have 
$$\frac{\int_{\Omega_\eps} |\sum_{i=1}^k c_i\grad I_\eps
u^{p_i}|^2}{\int_{\Omega_\eps} |\sum_{i=1}^k c_i I_\eps u^{p_i}|^2}
= \frac{ \langle A \sum_{i=1}^k c_i u^{p_i} | \sum_{i=1}^k c_i
u^{p_i}\rangle_{L^2}}{\|\sum_{i=1}^k c_i u^{p_i}\|_{L^2}^2}+o(1)~,$$
where the remainder $o(1)$ is uniform with respect to $c$ when $\eps$ goes to
zero. Using the Min-Max Principle another time, we get
$$ \min_{p_1<p_2<\ldots<p_k}~ \max_{c\in\RR^k_*}~
\frac{\int_{\Omega_\eps} |\sum_{i=1}^k c_i\grad I_\eps
u^{p_i}|^2}{\int_{\Omega_\eps} |\sum_{i=1}^k c_i I_\eps
u^{p_i}|^2}=\lambda^k+o(1)~.$$
This finishes the proof of Proposition \ref{prop1}.
%
\end{demo}


\section{Upper-semicontinuity of the spectrum}\label{up}

Let $\eps>0$ and let $(u^k_\eps)$ be a sequence of eigenfunctions of
$\Delta_\eps$ corresponding to the eigenvalues $(\lambda_\eps^k)$. We can assume
that the functions $u^k_\eps$ are orthogonal in $L^2(\Omega_\eps)$ and that
$\|u^k_\eps\|_{L^2}=\eps$. To work on a fixed domain, we set $\Omega=(0,1)^3$
and we introduce the functions $v^k_\eps=Ju^k_\eps$ where $J$ is the
canonical embedding of $H^1(\Omega_\eps)$ into $H^1(\Omega)$, that is that
$$Ju^k_\eps(y)=v^k_\eps(y)=v^k_\eps(y_1,\tilde y)=u^k_\eps(y_1,\eps \tilde
y)~.$$
We have 
$$ -\left(\partial^2_{y_1y_1}+\frac 1{\eps^2} \partial^2_{\tilde y\tilde
y}\right) v^k_\eps = \lambda^k_\eps v^k_\eps~~\text{ and }~~\|v^k_\eps
\|_{L^2}=1~.$$
By multiplying the previous equation by $v^k_\eps$ and integrating, we get
\begin{equation}\label{eq-v}
\int_\Omega |\partial_{y_1} v^k_\eps|^2 + \frac1{\eps^2}  |\partial_{\tilde y}
v^k_\eps|^2=\lambda^k_\eps~.
\end{equation}
Proposition \ref{prop1} shows that $(\lambda^k_\eps)_{\eps>0}$ is bounded.
Therefore, up to extracting a subsequence, we may assume that $(\lambda^k_\eps)$
converges to $\lambda^k_0=\liminf_{\eps\rightarrow 0}
\lambda^k_\eps$ when $\eps$ goes to $0$ and that
$(v^k_\eps)$ converges to a function $v^k_0\in H^1(\Omega)$, strongly in
$H^{3/4}(\Omega)$ and weakly in $H^1(\Omega)$. Moreover, \eqref{eq-v} shows that
$v^k_0$ depends only on $y_1$. In the following, we will abusively denote by
$v^k_0$, either the function in $H^1(\Omega)$ or the one-dimensional function in
$H^1(0,1)$.

The purpose of this section is to use the methods of \cite{Murat} (see also
\cite{Murat-CRAS1} and \cite{Murat-CRAS2}) to prove the following result.
\begin{prop}\label{prop2}
For all $k\in\NN^*$, the function $v^k_0$ is an eigenfunction of $A$ for the
eigenvalue $\lambda^k_0$.
\end{prop}

Proposition \ref{prop2} finishes the proof of Theorem \ref{th} since we
immediately get the upper-semicontinuity of the spectrum.
\begin{coro}
$$\forall k\in\NN^*~,~~\displaystyle \liminf_{\eps\rightarrow 0}
\lambda^k_\eps \geq \lambda^k~.$$
\end{coro}
\begin{demo}
We recall that the functions $v^k_\eps$ are orthonormalised in $L^2(\Omega)$
and converge strongly in $L^2(\Omega)$ to $v^k_0$. Thus, the functions $v^k_0$
are also orthonormalised. Since
$\lambda^k_0=\liminf_{\eps\rightarrow 0}
\lambda^k_\eps$, we know that
$\lambda^1_0\leq\lambda^2_0\leq\ldots\leq\lambda^k_0$. Then, Proposition
\ref{prop2} shows that $\lambda^1_0$, $\ldots$, $\lambda^k_0$ are $k$
eigenvalues of $A$ with linearly independent eigenfunctions, and thus that the
largest one $\lambda^k_0$ is larger than $\lambda^k$.
\end{demo}

The proof of Proposition \ref{prop2} splits into several lemmas. To simplify the
notations, we will omit the exponent $k$ in the remaining part of this
section and we will write $u_\eps$ for $u_\eps^k$, $v_\eps$ for $v_\eps^k$ etc.
\begin{lemma}\label{lemme2}
Let $B_\eps\subset\Omega_\eps$ be any cube of size $\eps$ and let $\Gamma_\eps$
be one of its faces. Then,
\begin{equation}\label{eq-lemme2}
\frac1{\eps^3} \int_{B_\eps} u_\eps(x)dx = \frac1{\eps^2} \int_{\Gamma_\eps}
u_\eps(\tilde x)d\tilde x+o(1)~.
\end{equation}
As a consequence, $v_0$ satisfies both Dirichlet boundary conditions
$v_0(0)=v_0(1)=0$.
\end{lemma}
\begin{demo}
We split the cube in slices $B_\eps=\cup_{s\in[0,\eps]} \Gamma_\eps(s)$ with
$\Gamma_\eps=\Gamma_\eps(0)$ and we set $x=(s,\tilde x)$ with $\tilde
x\in\Gamma_\eps(s)$. 
For each $s$, we have 
\begin{align*}
\left|\int_{\Gamma_\eps(s)}u_\eps(s,\tilde x)d\tilde x -
\int_{\Gamma_\eps(0)} u_\eps(0,\tilde x)d\tilde x \right| & \leq
\int_{\Gamma_\eps(\xi)} \int_0^s \left|\grad u_\eps(\xi,\tilde x)\right|
d\xi d\tilde x \\
&\leq \eps\sqrt{s} \sqrt{\int_{\Gamma_\eps(\xi)} \int_0^s \left|\grad
u_\eps(\xi,\tilde x)\right|^2
d\xi d\tilde x}~.
\end{align*}
To show \eqref{eq-lemme2}, we integrate the above inequality from $s=0$ to
$s=\eps$ and we notice that $\|\grad u_\eps\|_{L^2}=\lambda_\eps
\|u_\eps\|_{L^2}=\lambda_\eps \eps=\Oc(\eps)$.

The fact that $v_0(1)=0$ follows from $v_\eps(1,\tilde x)=0$ and the
strong convergence of $v_\eps$ to $v_0$ in $H^{3/4}(\Omega)$. To obtain the
other Dirichlet boundary condition, we apply \eqref{eq-lemme2} to the cube
$B_\eps=[0,\eps]\times [-\eps,0]\times [-\eps/2,\eps/2]$ at the left-end of
$\Omega_\eps$. Since $u_\eps$ vanishes on the upper face of $B_\eps$, the
average of $u_\eps$ goes to zero in $B_\eps$. Applying \eqref{eq-lemme2} again,
the average of $u_\eps$ goes to zero on the left face of $B_\eps$. Thus, the
average of $v_\eps$ goes to zero on the left face $\Gamma=\{0\}\times
[-1,0]\times [-1/2,1/2]$ of $\Omega$ and hence $\int_\Gamma v_0(0,\tilde
y)d\tilde y=0$ because $v_\eps$ converges to $v_0$ in $H^{3/4}(\Omega)$. Since
$v_0$ does not depend on $\tilde y$, this yields $v_0(0)=0$.
\end{demo}

We now focus on what happens close to the hole at $(a,0,0)$. To this end, we
use the notations of Section \ref{zoom} and we introduce the functions
$w_\eps\in H^1(K_\eps)$ defined by
$$\forall x \in K_\eps~,~~w_\eps(x)=u_\eps((a,0,0)+\delta\eps^2 x)~.$$
The functions $w_\eps$ will be useful to study the behaviour of
$u_\eps$ in the cube 
$$B_\eps=(a-\eps/2,a+\eps/2)\times(-\eps,0)\times(-\eps/2,
\eps/2)~.$$ 
We show that they weakly converges to $v_0(a)\zeta$ in $\dot H^1(K)$ in the
following sense.
\begin{lemma}\label{lemme3}
For all $\varphi\in \dot H^1(K)$, 
$$\int_{K_\eps} \grad w_\eps\grad\varphi \xrightarrow[~~\eps\longrightarrow
0~~]{} v_0(a) \int_K \grad \zeta \grad \varphi~.$$
\end{lemma}
\begin{demo}
We have  
$$\int_{K_\eps}|\grad w_\eps|^2=\frac1{\delta\eps^2}\int_{B_\eps}|\grad
 u_\eps|^2\leq \frac1{\delta\eps^2} \int_{\Omega_\eps}|\grad u_\eps|^2 =
 \frac1{\delta\eps^2} \lambda_\eps \int_{\Omega_\eps}|u_\eps|^2 =
 \frac{\lambda_\eps} {\delta}~.$$
Moreover, the average of $w_\eps$ in $K_\eps$ is equal to the one of
$u_\eps$ in $B_\eps$, which converges to $v_0(a)$ due to Lemma \ref{lemme2} and
the convergence of $v_\eps$ to $v_0$ in $H^{3/4}(\Omega)$. Applying Proposition
\ref{prop-weak}, we obtain the weak convergence of a subsequence of $w_\eps$ to
a limit $w_0$, whose average is $\overline w_0=v_0(a)$. To prove Lemma
\ref{lemme3}, it remains to show that $w_0=v_0(a)\zeta$, which does not depend
on the chosen subsequence $(\eps_n)$.

Let $\varphi\in\Cc^\infty_0(K)$ and assume that $\eps$ is small enough such
that $\supp(\varphi)\subset K_\eps$. We set  
$$\forall x\in
B_\eps~,~~\tilde\varphi_\eps(x)=\varphi\left(\frac{x-(a,0,0)}{\delta\eps^2}
\right)$$
and we extend $\tilde\varphi_\eps$ by zero in $\Omega_\eps$. Since 
$$\| u_\eps \|_{L^2}=\eps~\text{ and
}~\|\tilde\varphi_\eps\|_{L^2(\Omega_\eps)}=\delta^{3/2} \eps^3
\|\varphi\|_{L^2(K)}~,$$
we get
$$\int_{K_\eps} \grad w_\eps\grad\varphi= \frac 1{\delta \eps^2} \int_{B_\eps}
\grad u_\eps \grad \tilde\varphi_\eps= \frac 1{\delta \eps^2}
\int_{\Omega_\eps} \Delta u_\eps \tilde\varphi_\eps= \frac
{\lambda_\eps}{\delta \eps^2} \int_{\Omega_\eps} u_\eps \tilde\varphi_\eps
\xrightarrow[~~\eps\longrightarrow
0~~]{} 0~.$$
Thus, $w_0$ is orthogonal to $\Cc^\infty_0(K)$ and hence to $\dot
H^1_0(K)$ and Proposition \ref{prop-redef-scalar} implies that $w_0=\overline
w_0\zeta$. Since we already know that $\overline w_0=v_0(a)$, Lemma \ref{lemme3}
is proved.
\end{demo}

{ \noindent \emph{\textbf{Proof of Proposition \ref{prop2}~:}}}
We have shown in Lemma \ref{lemme2} that $v_0$ satisfies Dirichlet boundary
condition at $x_1=0$ and $x_1=1$. Let $\phi\in H^1_0(0,1)$ be a test function.
We also denote by $\phi$ the canonical embedding of $\phi$ into $H^1(\Omega)$. 
We embed $\phi$ into $\Omega_\eps$ by setting $\phi_\eps=I_\eps \phi$, where
$I_\eps$ is the embedding introduced in the proof of Proposition \ref{prop1}.
Using the notations of Figure \ref{Ieps}, we have
\begin{equation}\label{eq-fin}
\int_{\Omega_\eps} \grad u_\eps \grad \phi_\eps = \int_{x_1<a-\eps/2} \grad
u_\eps \grad \varphi_\eps + \phi(a)\int_{B_\eps} \grad u_\eps\grad \tilde
\zeta_\eps +  \int_{x_1>a+\eps/2} \grad
u_\eps \grad \psi_\eps
\end{equation}
The limits of the different terms are as follows. 
First, notice that 
$$\int_{\Omega_\eps} \grad u_\eps \grad \phi_\eps =
\lambda_\eps \int_{\Omega_\eps} u_\eps \phi_\eps = \eps^2\lambda_\eps
\int_{\Omega} v_\eps J\phi_\eps$$
where $J\phi_\eps$ is the canonical embedding
of $\phi_\eps$ in $H^1(\Omega)$. Obviously, $J\phi_\eps$ converges to $J\phi$ in
$L^2(\Omega)$ and we know that $v_\eps$ converges to $v_0$ in $L^2(\Omega)$.
Thus,
$$\int_{\Omega_\eps} \grad u_\eps \grad \phi_\eps
=\eps^2 \lambda_0 \int_\Omega v_0
\phi + o(\eps^2) =\eps^2 \lambda_0\int_0^1 v_0 \phi+o(\eps^2)~.$$
In the parts $x_1<a-\eps/2$ and $x_1>a+\eps/2$, we know that $v_\eps$ converges
to $v_0$ weakly in $H^1(\Omega)$ and obviously $J\varphi_\eps$ and $J\psi_\eps$
converge to $\phi$ strongly in $H^1$. Moreover, notice that $J\varphi_\eps$ and
$J\psi_\eps$ only depends on $x_1$. Hence, 
\begin{align*}
\int_{x_1<a-\eps/2} \grad
u_\eps \grad \varphi_\eps & + \int_{x_1>a+\eps/2} \grad
u_\eps \grad \psi_\eps\\ 
& = \eps^2\left( \int_{x_1<a-\eps/2} \partial_{x_1}
v_\eps \partial_{x_1} (J\varphi_\eps) + \int_{x_1>a+\eps/2} \partial_{x_1}
v_\eps \partial_{x_1} (J\psi_\eps) \right)\\
& = \eps^2 \left(
\int_0^{a} \partial_{x_1} v_0\partial_{x_1} \phi + \int_{a}^1
\partial_{x_1} v_0\partial_{x_1} \phi\right)+o(\eps^2)~.\\
& = \eps^2 \int_0^1 \partial_{x_1} v_0\partial_{x_1} \phi +o(\eps^2)~.
\end{align*}
The term of \eqref{eq-fin} in the box $B_\eps$ is more delicate, but all the
work has already been done in Lemma \ref{lemme3}. Indeed we have
$$\int_{B_\eps} \grad u_\eps\grad \tilde
\zeta_\eps=\delta \eps^2 \int_{K_\eps} \grad w_\eps \grad \zeta_\eps~.$$
By definition $\zeta_\eps$ converges to $\zeta$ strongly in $\dot H^1(K)$.
Thus, Lemma \ref{lemme3} implies that 
$$\int_{B_\eps} \grad u_\eps\grad \tilde
\zeta_\eps=\delta \eps^2 v_0(a) \int_{K} \grad \zeta \grad \zeta +o(\eps^2)=
\alpha\delta v_0(a) \eps^2 + o(\eps^2)~.$$
In conclusion, when $\eps$ goes to $0$, Equality \eqref{eq-fin} shows that 
$$ \lambda_0\int_0^1 v_0 \phi=\int_0^1 \partial_{x_1} v_0\partial_{x_1} \phi +
\alpha\delta v_0(a)\phi(a)~.$$
Since this holds for all $\phi\in H^1_0(0,1)$, going back to the variational
form of $A$ given in \eqref{eq-var-A}, this shows that $v_0$ is an
eigenfunction of $A$ for the eigenvalue $\lambda_0$ (remember that
$\|v_0\|_{L^2}=1$ and so $v_0$ is not zero).
{\hfill$\square$\\   \vspace{0.4cm}}



\begin{thebibliography}{99}

\bibitem{Anne}
C. Ann\'e, {\it Spectre du laplacien et \'ecrasement d'anses},
Annales Scientifiques de l'\'Ecole Normale Sup\'erieure \no 20 (1987),
pp. 271-280. 

\bibitem{AHH}
J.M. Arrieta, J.K. Hale and Q. Han, {\it 
Eigenvalue problems for non-smoothly perturbed domains},
Journal of Differential Equations \no 91 (1991), pp. 24-52. 

\bibitem{Beale}
J.T. Beale, {\it Scattering frequencies of resonators}, Communications on Pure
and Applied Mathematics \no 26 (1973), pp. 549-563.

\bibitem{Benade}
A.H. Benade, {\it On the Mathematical Theory of Woodwind Finger Holes},
Journal of the Acoustical Society of America, \no 32 (1960), pp.
1591-1608.

\bibitem{Bolton}
P. Bolton, {\rm\tt http://www.flute-a-bec.com/acoustiquegb.html}, the website of
a recorder maker.

\bibitem{Ciuperca}
I.S. Ciuperca, {\it Reaction-diffusion equations on thin domains with varying
order of thinness}, Journal of Differential Equations \no 126 (1996),
pp. 244-291. 

\bibitem{Coltman}
J.W. Coltman, {\it Acoustical analysis of the Boehm flute}, Journal of the
Acoustical Society of America, \no 65 (1979), pp. 499-506. 

\bibitem{CH}
R. Courant and D. Hilbert, {\it Methods of mathematical physics. Vol. I}.
Interscience Publishers, New York, 1953. 

\bibitem{Dickens}
P.A. Dickens, {\it Flute acoustics: measurements, modelling and design}. PhD
Thesis, University of New South Wales, 2007.

\bibitem{Evans}
L.C. Evans, {\it Partial differential equations}. Graduate Studies in
Mathematics \no 19. American Mathematical Society, Providence, RI, 1998.

\bibitem{Fletcher}
N.H. Fletcher and T.D. Rossing, {\it The Physics of Musical
Instruments}. Springer-Verlag, New York, 1998. 

\bibitem{Hadamard}
J. Hadamard, {\it La th\'eorie des plaques \'elastiques planes},
Transactions of the American Mathematical Society \no 3 (1902), pp. 401-422.

\bibitem{Hale-Raugel}
J.K. Hale and G. Raugel, {\it Reaction-diffusion equation on thin
domains}, Journal de Math\'ematiques Pures et Appliqu\'ees \no 71 (1992), pp.
33-95. 

\bibitem{JM}
S. Jimbo and Y. Morita, {\it Remarks on the behavior of certain
eigenvalues on a singularly perturbed domain with several thin channels},
Communications in Partial Differential Equations \no 17 (1992), pp.
523-552.

\bibitem{LS}
M. Lobo and E. S\'anchez-Palencia, {\it Sur certaines propri\'et\'es spectrales
des perturbations du domaine dans les probl\`emes aux limites},
Communication in Partial Differential Equations \no 4 (1979), pp.
1085-1098. 

\bibitem{Murat-CRAS1}
J. Casado-D\'\i az, M. Luna-Laynez and F. Murat, {\it
Asymptotic behavior of diffusion problems in a domain made of two cylinders of
different diameters and lengths},
Comptes Rendus Math\'ematique. Acad\'emie des Sciences. Paris
\no 338 (2004), pp. 133-138. 

\bibitem{Murat-CRAS2}
J. Casado-D\'\i az, M. Luna-Laynez and F. Murat, {\it
Asymptotic behavior of an elastic beam fixed on a small part of one of its
extremities},
Comptes Rendus Math\'ematique. Acad\'emie des Sciences. Paris 
\no 338 (2004), pp. 975-980. 


\bibitem{Murat}
J. Casado-D\'\i az, M. Luna-Laynez and F. Murat, {\it
The diffusion equation in a notched beam},
Calculus of Variations and Partial Differential Equations \no 31 (2008),
pp. 297-323. 

\bibitem{PR}
M. Prizzi and K. Rybakowski, {\it The effect of domain squeezing upon the
dynamics of reaction-diffusion equations}, Journal of Differential Equations
\no 173 (2001), pp. 271-320. 
 
\bibitem{Raugel}
G. Raugel, {\it Dynamics of partial differential equations on thin domains}.
Dynamical systems (Montecatini Terme, 1994). Lecture Notes in
Mathematics \no 1609, pp. 208-315 . Springer, Berlin, 1995. 
 
\bibitem{Reed-Simon}
M. Reed and B. Simon, {\it Methods of Modern Mathematical Physics IV: Analysis
of Operators}. Academic Press, 1978.

\bibitem{Rossing}
T.D. Rossing, {\it The Science of Sound}. Addison-Wesley,
Reading, Mass, 1982.

\bibitem{Schatzman}
M. Schatzman, {\it On the eigenvalues of the Laplace operator on a thin set
with Neumann boundary conditions}, Applicable Analysis \no 61 (1996), pp.
293-306.

\bibitem{Simader-Sohr}
C.G. Simader and H. Sohr, {\it The Dirichlet problem for the Laplacian in
bounded and unbounded domains}. Pitman Research Notes in Mathematics Series \no
360. Longman, Harlow, 1996.

\bibitem{Wolfe-web}
J. Wolfe, {\rm\tt http://www.phys.unsw.edu.au/jw/fluteacoustics.html}, the
website of an acoustician.

\bibitem{Wolfe}
J. Wolfe and J. Smith, {\it Cutoff frequencies and cross fingerings in baroque,
classical, and modern flutes}, Journal of the Acoustical Society of America \no
114 (2003), pp. 2263-2272.

\end{thebibliography}
\end{document}